\documentclass[12pt]{article}

\usepackage{amsmath}
\usepackage{amssymb}
\usepackage{amscd}

\newtheorem{thm}{Theorem}[section]

\newtheorem{lem}[thm]{Lemma}

\newcommand{\gp}{\mathfrak{p}}
\newcommand{\gq}{\mathfrak{q}}

\def\ZZ{\mathbb{Z}}

\newcommand{\gr}[1]{{\mathrm G}(#1)}

\def\deg{{\rm deg \,}}

\def\ker{{\rm Ker}}

\def\rank{{\rm rank}}

\def\dep{{\rm depth \,}}
\def\gr{{\rm grade \,}}

\newcommand{\ass}[2]{{\rm Ass}_{#1}\,#2}

\newcommand{\som}[3]{{\rm M}(#1, #2\,; #3)} 
\newcommand{\deti}[2]{{\rm I}_{#1}(#2)} 

\begin{document}

\title{\Large On the symmetric and Rees algebras \\ of certain determinantal ideals \\}

\author{Kosuke Fukumuro}
\date{
\small
Graduate School of Science, Chiba University, \\
1-33 Yayoi-Cho, Inage-Ku, Chiba-Shi, 263-8522, JAPAN}

\maketitle

\section{Introduction}
Let $R$ be a Noetherian ring and let $m$, $n$ be integers such that $1 \leq m \leq n$. 
We denote by
$\som{m}{n}{R}$ the set of $m \times n$ matrices with entries in $R$. 
Let $M$ = $(x_{ij}) \in \som{m}{n}{R}$ and
$S$ a polynomial ring over $R$ with variables $T_1,T_2, \dots , T_n$. 
We regard $S$ as a graded ring
by setting $\deg$ $T_j = 1$ for  $1 \leq  \forall j \leq n$. Let
\[
f_i = \displaystyle \sum_{j = 1}^n x_{ij}T_j \in S_1
\]
for $1 \leq \forall i \leq m$ 
and let $\deti{k}{M}$ be the ideal of $R$ generated by the $k$-minors of $M$
for  $1 \leq \forall k \leq m$.  
The purpose of this paper is to give elementary proofs to the following
theorems due to Avramov \cite{a}.

\begin{thm}\label{1.1}
{\rm \cite{a}}. The following conditions are equivalent:
\begin{itemize}
\item[{\rm (1)}]$\gr \deti{k}{M} \geq m - k +1$ \, for $1 \leq \forall k \leq m$.
\item[{\rm (2)}]$\gr (f_1,f_2, \dots ,f_m)S = m$.
\end{itemize}
\end{thm}

\begin{thm}\label{1.2}
{\rm \cite{a}}. 
Suppose $n = m + 1$. We set $I = \deti{m}{M}$. Let ${\rm S}(I)$ and ${\rm R}(I)$ be the
symmetric algebra and the Rees algebra of $I$, respectively. 
Then the following conditions are equivalent:
\begin{itemize}
\item[{\rm (1)}]$\gr \deti{k}{M} \geq m - k + 2$ \,\, for $1 \leq \forall k \leq m$.
\item[{\rm (2)}]
{\rm (i)} The natural map ${\rm S}(I) \to {\rm R}(I)$ is isomorphic,  
{\rm (ii)} ${\rm S}(I) \cong S/(f_1, f_2, \dots ,f_m)S$ as  graded $R$-algebras and
{\rm (iii) }$\gr (f_1 , f_2 , \dots , f_m)S = m$\,.
\end{itemize}

\end{thm}

The Rees algebra of $I$ is the subalgebra of the polynomial ring 
$R[t]$ generated by $It$ over R.
If the conditions (i) and (ii) of 1.2 are satisfied, there exists a surjection
$S \to {\rm R}(I)$ of graded $R$-algebras 
whose kernel coincides with $(f_1,f_2,\dots,f_m)S$.
Moreover, if the condition (iii) of 1.2 is
satisfied, the Koszul complex of $f_1,f_2, \dots ,f_m$ 
is an acyclic complex of graded free $S$-modules
by \cite[1.\,6.\,17]{bh}\,.
Therefore, under the condition (2) of 1.2, 
we get a graded $S$-free resolution of ${\rm R}(I)$, and
taking its homogeneous part of degree $r \in \ZZ$, 
we get an $R$-free resolution of $I^r$, from which we
can deduce some homological properties of powers of $I$.
In the subsequent paper \cite{fin}, using the $R$-free resolution of $I^r$ 
constructed in this way,
we study the associated prime ideals of $R/{I^r}$ and
compute the saturation of $I^m$.
So, the author thinks that 1.2 is very convenient and it may have more
application.
Although we may lose sight of some important meaning of the original proof,
however, the existence of an elementary proof is quite helpful for users.

\section{Preliminaries}
In this section we summarize preliminary results. Although these facts might be well-known, we give the proofs for completeness.

\begin{lem}\label{2.1}
$\ass{}{S} = \{\, \gp S \,\mid\, \gp \in \ass{}{R} \,\}$.
\end{lem}

\noindent
{\it Proof}.\,
It is enough to show in the case where $n = 1$. We put $T = T_1$.

Let us take any $\gp \in \ass{}{R}$.
Then there exists $x \in R$ such that $\gp = 0 :_R x$.
It is easy to see
$\gp S = 0 :_S x$. On the other hand, $S/{\gp S} \cong (R/\gp)[T]$, 
which is an integral domain. Hence
$\gp S \in \ass{}{S}$.

Conversely, we take any $Q \in \ass{}{S}$.
Then $Q$ is homogeneous ,
and so $Q = 0 :_S yT^n$ for some $y \in R$ and $0 \leq n \in \ZZ$.
We put $\gq = Q \cap R$.
Then $\gq \in {\rm Spec}\,R$ and $\gq = 0 :_R y$,
which means $\gq \in \ass{}{R}$.
Moreover, we have $\gq S = 0 :_S yT^n = Q$.
Thus the proof is complete.

\begin{lem}\label{2.2}
$f_1$ is a non-zerodivisor on $S$ if and only if 
$\gr (x_{11},x_{12}, \dots , x_{1n})R > 0$.
\end{lem}

\noindent
{\it Proof}.\,
Suppose $\gr (x_{11},x_{12}, \dots , x_{1n})R = 0$.
Then there exists $\gp \in \ass{}{R}$ such that
$(x_{11},x_{12}, \dots , x_{1n})R \subseteq \gp$. In this case, 
we have $f_1\in \gp S \in \ass{}{S}$, and so
$f_1$ is a zerodivisor on $S$.

Conversely, suppose that $f_1$ is a zerodivisor on $S$.
Then there exists $Q \in \ass {}{S}$ such that
$f_1 \in Q$, and $Q = \gq S$ for some $\gq \in \ass{}{R}$.
In this case, we have $x_{1j} \in \gq$ for
$1 \leq \forall j \leq n$, which means $\gr (x_{11},x_{12}, \dots , x_{1n})R = 0$. 
Thus the proof is complete.

\begin{lem}\label{2.3}
Let $N = (y_{ij}) \in \som{m}{n}{R}$ be the matrix 
induced from $M$ by some elementary operation. We set
\[
g_i = \displaystyle \sum_{j=1}^n y_{ij}T_j
\]
for  $1 \leq \forall i \leq m$. 
Then there exists an isomorphism 
$\varphi : S \stackrel{\sim}{\to} S$ of $R$-algebras such that 
$\varphi ((f_1,f_2, \dots , f_m)S) = (g_1,g_2, \dots , g_m)S$. 
In particular, we have
$\gr (f_1,f_2, \dots , f_m)S = \gr (g_1,g_2, \dots , g_m)S$.
\end{lem}

\noindent
{\it Proof}. \,
The elementary operation is one of the followings:
\begin{itemize}
\item[{\rm (i)}] Fixing a unit $c \in R$ and $1 \leq i \leq m$, 
replace $x_{ij}$ with $cx_{ij}$ for $1 \leq \forall j \leq n$.
\item[{\rm (ii)}] Fixing an element $c \in R$ and $1 \leq i, k \leq m \, (i \neq k)$, 
replace $x_{ij}$ with $x_{ij} + cx_{kj}$ for  $1 \leq \forall j \leq n$.
\item[{\rm (iii)}] Fixing $1 \leq i < k \leq m$, exchange 
$x_{ij}$ with $x_{kj}$ for  $1 \leq \forall j \leq n$.
\item[{\rm (iv)}] Fixing a unit $c \in R$ and $1 \leq j \leq n$, 
replace $x_{ij}$ with $cx_{ij}$ for $1 \leq \forall i \leq m$.
\item[{\rm (v)}] Fixing an element $c \in R$ and $1 \leq j$, 
$\ell \leq n \, (j \neq \ell)$, replace $x_{ij}$ 
with $x_{ij} + cx_{i \ell}$ for $1 \leq \forall i \leq m$.
\item[{\rm (vi)}] Fixing $1 \leq j < \ell \leq n$, exchange 
$x_{ij}$ with $x_{i \ell}$ for $1 \leq \forall i \leq m$.
\end{itemize}
\noindent
In each case stated above the following relations between $f_i '$s and $g_i '$s hold: 
\begin{itemize}
\item[{\rm (i)}]$g_i = cf_i$ and $g_p = f_p$ for any $p \neq i$.
\item[{\rm (ii)}]$g_i = f_i + cf_k$ and $g_p = f_p$ for any $p \neq i$.
\item[{\rm (iii)}]$g_i = f_k$, $g_k = f_i$ and $g_p = f_p$ for any $p \neq i, k$.
\item[{\rm (iv)}]$g_i = f_i (T_1, \dots, T_{j-1}, cT_j, T_{j+1}, \dots, T_n)$ for any $i$.
\item[{\rm (v)}]$g_i = f_i (T_1, \dots, T_{j-1}, T_j +cT_\ell, T_{j+1}, \dots, T_n)$ for any $i$.
\item[{\rm (vi)}]$g_i = f_i(T_1, \dots ,\overset{\overset{j}{\smallsmile}}{T_\ell}, \dots ,
\overset{\overset{\ell}{\smallsmile}}{T_j}, \dots , T_n)$
\end{itemize}
\noindent
In the cases of (i), (ii) and (iii), we set $\varphi = {\rm id}_S$. In the other cases we define $\varphi$ by
\begin{itemize}

\item[{\rm (iv)}]$\varphi (T_j) = cT_j$ and  $\varphi (T_q) = T_q$ for any $q \neq j$,
\item[{\rm (v)}]$\varphi (T_j) = T_j + cT_\ell$ and $\varphi (T_q) = T_q$ for any $q \neq j$,
\item[{\rm (vi)}]$\varphi (T_j) = T_\ell,  \varphi (T_\ell) = T_j$ and $ \varphi (T_q) = T_q$ for any 
$q \neq j, \ell$.

\end{itemize}
Then, $\varphi$ becomes an isomorphism and the required equality holds.

\vspace{1em}

Let $K_{\bullet}$ be the graded Koszul complex with respect to $f_1, f_2, \dots, f_m$. 
By $\partial_\bullet$ we denote its boundary map . 
Let $e_1, e_2, \dots, e_m$ be an $S$-free basis of $K_1$ 
consisting of homogeneous elements of degree 1 
such that $\partial_1(e_i) = f_i$ for  $1 \leq \forall i \leq m$. 
Then, for $1 \leq \forall r \leq m$,
\[
{\{e_{i_1}\wedge \cdots \wedge e_{i_r} | 1 \leq i_1 < \cdots < i_r \leq m} \}
\]
is an $S$-free basis of $K_r$ 
consisting of homogeneous elements of degree $r$, and we have
\[
\partial_r (e_{i_1}\wedge \cdots \wedge e_{i_r}) =
\sum_{p=1}^{r} (-1)^{p-1} f_{ip} \cdot e_{i_1} \wedge \cdots \wedge \widehat{e_{i_p}} \wedge \cdots \wedge e_{i_r}.
\]
Let $\ell \in \ZZ$. 
Taking the homogeneous part of degree $\ell$ of $K_{\bullet}$, 
we get a complex
\[
[K_{\bullet}]_{\ell} : 0 \longrightarrow [K_m]_{\ell} \stackrel{[\partial_m]_{\ell}}{\longrightarrow}
[K_{m-1}]_{\ell} \longrightarrow \cdots \longrightarrow [K_1]_{\ell} \stackrel{[\partial_1]_{\ell}}{\longrightarrow} [K_0]_{\ell} \longrightarrow 0
\]
of finitely generated free $R$-modules, 
where $[\partial_r]_{\ell}$ denotes the restriction of $\partial_r$ to
$[K_r]_{\ell}$ for $1 \leq \forall r \leq m$. 
It is obvious that $[K_r]_{\ell} = 0$ if $\ell < r$. 
On the other hand, if $\ell \geq r$, then
\begin{multline*}
\{\,
T_1^{\alpha_1} T_2^{\alpha_2} \cdots T_n^{\alpha_n} e_{i_1} \wedge \cdots \wedge e_{i_r}
\,
\mid
\,
0 \leq \alpha_1,\alpha_2, \cdots, \alpha_n \in \ZZ\,, \\
\sum_{k = 1}^n \alpha_k = l - r\,,
1\leq i_1 < i_2 < \cdots < i_r \leq m
\,\}
\end{multline*}
is an $R$-free basis of $[K_r]_{\ell}$.

\begin{lem}\label{2.4}
$\rank_R [K_m]_m = 1$, $\rank_R [K_{m-1}]_m = mn$ and $\deti{1}{[\partial_m]_m} = \deti{1}{M}.$
\end{lem}

\noindent
{\it Proof}. \,
We get $\rank_R [K_m]_m = 1$ as $[K_m]_m$ is generated by 
$e_1 \wedge e_2 \wedge \cdots \wedge e_m$. 
Moreover, we get $\rank_R [K_{m-1}]_m = mn$ 
as $\{{T_j \check{e}_i \, | \, 1 \leq i \leq m, 1 \leq j \leq n}\}$ 
is an $R$-free basis of $[K_{m-1}]_m$,
where $\check{e}_i = e_1 \wedge \cdots \wedge \widehat{e_i} \wedge \cdots \wedge e_m$ 
for $1 \leq \forall i \leq m$. 
The last assertion holds since

\begin{equation}
\begin{split}
\partial_m (e_1\wedge e_2 \wedge \cdots \wedge e_m) &= \displaystyle \sum_{i = 1}^m (-1)^{i-1} f_i
\check{e}_i \notag \\
& = \displaystyle \sum_{i=1}^m (-1)^{i-1} (\sum_{j=1}^n x_{ij}T_j)
\check{e}_i \notag \\
& = \displaystyle \sum_{i=1}^m \sum_{j=1}^n (-1)^{i-1}x_{ij} T_j \check{e}_i \notag. \\
\end{split}
\end{equation}

\begin{lem}\label{2.5}
$\rank_R [K_1]_1 = m$, $\rank_R [K_0]_1 = n$ and $\deti{m}{[\partial_1]_1} = \deti{m}{M}$.
\end{lem}

\noindent
{\it Proof}. \,
We get $\rank_R[K_1]_1 = m$ as $e_1, e_2, \cdots, e_m$ is an $R$-free basis of $[K_1]_1$. 
Moreover,
we get $\rank_R [K_0]_1 = n$ as $T_1, T_2, \cdots, T_n$ 
is an $R$-free basis of $[K_0]_1$. 
The last assertion holds since
\[
\partial_1 (e_i) = f_i = \displaystyle \sum_{j=1}^n x_{ij}T_j
\]
for $1 \leq \forall i \leq m$.

\section{Proofs of 1.1 and 1.2}

\noindent
{\it Proof of 1.1}. \,
We prove by induction on $m$. If $m=1$, then the assertion holds by 2.2.
So, we may assume $m \geq 2$. \\
 
$(1) \Rightarrow(2)$ (The proof of this implication is same as that of {\rm \cite[Propsosition 1]{a}}. 
Put $t = \gr(f_1, f_2, \dots, f_m)S$ and suppose $t < m$. 
Take a maximal $S$-regular sequence $g_1, g_2, \dots, g_t$ 
contained in $(f_1, f_2, \dots, f_m)S$. Then, 
as any element in $(f_1, f_2, \dots, f_m)S$ is a zerodivisor on 
$S/(g_1, g_2, \dots, g_t)S$, there exists $P \in \ass{S}{S/(g_1, g_2, \dots, g_t)S}$ 
such that $(f_1, f_2, \dots, f_m)S \subseteq P$. 
When this is the case, we have $\dep S_P = \gr P = t < m$.
Here we put $\gp = P \cap R$. 
Then $\gr \gp \leq \gr P < m$, and so $\deti{1}{M} \nsubseteq \gp$ as
$\gr \deti{1}{M} \geq m$. 
Hence, replacing rows and columns of $M$ if necessary, 
we may assume
$x_{11} \notin \gp$. 
Then, applying elementary operations to $M$ in $\som{m}{n}{R_{\gp}}$, 
we get a matrix of the form
\[
\left(\begin{array}{@{\,}c|ccc@{\,}}
1 & 0 & \cdots & 0\\ \hline
0 &&& \\
\vdots & &N &\\
0&&&
\end{array}
\right),
\]
where $N = (y_{ij}) \in \som{m-1}{n-1}{R_{\gp}}$. For $1 \leq \forall k \leq m-1$, 
we have $\deti{k}{N} = \deti{k+1}{M}R_{\gp}$, and so
\[
\gr \deti{k}{N} \geq \gr \deti{k+1}{M} \geq m-(k+1)+1 =(m-1)-k+1.
\]
For $2 \leq \forall i \leq m$, we put
\[
g_i = \displaystyle \sum_{j=2}^{n} y_{i-1,j-1} T_{j} \in S' ,
\]
where $S' = R_{\gp}[T_2, \cdots, T_n]$. 
By the hypothesis of induction, 
we have $\gr (g_2, \dots, g_m)S' = m-1$. Then, 
as $S_{\gp} = S'[T_1]$,
it follows that $\gr (T_1, g_2, \dots, g_m)S_{\gp} = m$. 
Hence we get $\gr (f_1, f_2, \dots, f_m)S_{\gp} = m$ by 2.3, 
and so $\gr (f_1, f_2, \dots, f_m)S_P = m$, 
which contradicts to $\dep S_P < m$. 
Thus we see $t = m$. \\

$(2) \Rightarrow (1)$ By\cite[1.6.17]{bh}, 
$K_{\bullet} = K_{\bullet}(f_1, \dots, f_m)$ is acyclic, and so
\[
0 \longrightarrow [K_m]_m \stackrel{[\partial_m]_m}{\longrightarrow}
[K_{m-1}]_m \longrightarrow \cdots \longrightarrow [K_1]_m \stackrel{[\partial_1]_m}{\longrightarrow} [K_0]_m \longrightarrow 0
\]
is acyclic, too. 
Then, as $[K_m]_m \cong R$ and $\deti{1}{[\partial_m]_m} = \deti{1}{M}$ 
by 2.4, we get $\gr \deti{1}{M} \geq m$ by \cite[1.4.13]{bh}. 
Suppose that $\ell := \gr \deti{k}{M} \leq m-k$ 
for some $k$ with $2 \leq k \leq m$.
We take a maximal $R$-regular sequence $c_1, c_2, \dots, c_{\ell}$ 
contained in $\deti{k}{M}$.
Then, as any element in $\deti{k}{M}$ is a zerodivisor on $R/(c_1, c_2, \dots, c_\ell)R$, 
we have $\deti{k}{M} \subseteq \gp$ for some $\gp \in \ass{R}{R/(c_1, c_2, \dots, c_\ell)R}$. 
When this is the case, we
have $\dep R_{\gp} = \gr \gp = \ell \leq m-k$, which means $\deti{1}{M} \nsubseteq \gp$. 
Then, applying elementary operations to $M$ in $\som{m}{n}{R_{\gp}}$,
 we get a matrix of the form
\[
\left(\begin{array}{@{\,}c|ccc@{\,}}
1 & 0 & \cdots & 0\\ \hline
0 &&& \\
\vdots & &N &\\
0&&&
\end{array}
\right),
\]
where $N = (y_{ij}) \in \som{m-1}{n-1}{R_{\gp}}$. 
For $2 \leq \forall i \leq m$, 
we put
\[
g_i = \displaystyle \sum_{j=2}^{n} y_{i-1,j-1} T_{j} \in S' ,
\]
where $S' = R_{\gp}[T_2, \cdots, T_n]$. 
By 2.3, we have 
\[
\gr (T_1, g_2, \dots, g_m)S_{\gp} = \gr (f_1, f_2, \dots, f_m)S_{\gp} = m.
\]
Hence $\gr (g_2, \dots, g_m)S' = m - 1$ since $S_{\gp} = S'[T_1]$.
Then, the hypothesis of induction implies
\[
\gr \deti{k-1}{N} \geq (m-1)-(k-1)+1 = m-k+1,
\]
which contradicts to $\dep R_{\gp} \leq m-k$ as $\deti{k-1}{N} = \deti{k}{M}R_{\gp}$.
Thus we get 
\[
\gr \deti{k}{M} \geq m-k+1 
\] 
for $1 \leq \forall k \leq m$ and the proof is complete.

\hspace{1em}

\noindent
{\it Proof of 1.2}. \,
$(1) \Rightarrow (2)$\label{1-2}
By the hypothesis we have $\gr I \geq 2$, and so
\[
0 \longrightarrow R^m \stackrel{{}^t \! M}{\longrightarrow} R^{m+1} \longrightarrow I \longrightarrow 0
\]
is a free resolution of $I$ by \cite[1.4.17]{bh}. Hence the condition (ii) is satisfied.
Moreover, we see that the condition (iii) is satisfied 
by $(1) \Rightarrow (2)$ of 1.1. 
So, we have to show the assertion of (i).
 Let $L$ be the kernel of the natural map $S(I) \to R(I)$, 
 and consider the exact sequence
\[
0 \longrightarrow L \longrightarrow {\rm S}(I) \longrightarrow {\rm R}(I) \longrightarrow 0
\]
of $S$-modules. 
Let us prove $L = 0$ by induction on $m$.
  
First, we consider the case where $m=1$. Suppose $L \neq 0$.
Then there exists $P \in \ass{S}{L}$. 
Because $L \subseteq {\rm S}(I) \cong S/f_1 S$ by (ii),
we have $P \in \ass{S}{S/f_1 S}$, and so $\gr P = 1$ by (iii).
We put $\gp = P \cap R$. Then $\gr \gp \leq 1$. 
Because $I = \deti{1}{M}$ and $\gr \deti{1}{M} \geq 2$,
we have $I \nsubseteq \gp$. 
This means $I_{\gp} = R_{\gp}$, 
and so the natural map ${\rm S}(I_{\gp}) \to {\rm R}(I_{\gp})$ is an isomorphism. 
Then, looking at the commutative diagram

$$
\begin{CD}
0 @>>> L_{\gp} @>>> {\rm S}(I)_{\gp} @>>> {\rm R}(I)_{\gp} @>>> 0 \,\, {\rm (ex)} \\
@. @. @VV{\wr}V @VV{\wr}V \\
@.@. {\rm S}(I_{\gp}) @>{\sim}>> {\rm R}(I_{\gp}) ,
\end{CD}
$$
we get $L_{\gp} = 0$, and so $L_P = 0$, which contradicts to $P \in \ass{R}{L}$. Therefore we see $L =
0$.

Next, we consider the case where $m \geq 2$. 
Suppose $L \neq 0$. Then there exists $P \in \ass{S}{L}$.
Because $L \subseteq {\rm S}(I) \cong S/(f_1, \dots, f_m)S$ by (ii), 
we have $P \in \ass{S}{S/(f_1, \dots, f_m)S}$, and so $\gr P = m$ by (iii). 
We put $\gp = P \cap R$.
Then $\gr \gp \leq m$, and so $\deti{1}{M} \nsubseteq \gp$ as $\gr \deti{1}{M} \geq m+1$.
Hence, applying elementary operations to $M$ in $\som{m}{m+1}{R_{\gp}}$, 
we get a matrix of the form
\[
\left(\begin{array}{@{\,}c|ccc@{\,}}
1 & 0 & \cdots & 0\\ \hline
0 &&& \\
\vdots & &N &\\
0&&&
\end{array}
\right),
\]
where $N = (y_{ij}) \in \som{m-1}{m}{R_{\gp}}$. 
When this is the case, for  $1 \leq \forall k \leq m-1$, 
we have

\begin{equation}
\begin{split}
\gr \deti{k}{N} & = \gr \deti{k+1}{M} R_{\gp} \\
& \geq m-(k+1)+2 \\
& = (m-1)-k+2 \,\,\, . \notag \\
\end{split}
\end{equation}
We notice that $\deti{m-1}{N} = \deti{m}{M}R_{\gp} = I_{\gp}$, 
and so the hypothesis of induction implies that
the natural map ${\rm S}(I_{\gp}) \to {\rm R}(I_{\gp})$ is isomorphic. 
Now we look at the commutative diagram above again, 
and get $L_{\gp} = 0$. 
Hence $L_P = 0$, which contradicts to $P \in \ass{S}{L}$. \\

$(2) \Rightarrow (1)$ \label{2-1}
By (i) and (ii), there exists a surjection $\pi : S \to {\rm R}(I)$ 
of graded $R$-algebras such that
$\ker \, \pi = (f_1, f_2, \dots, f_m)S$. 
On the other hand, $K_{\bullet} = K_{\bullet} (f_1, f_2, \dots, f_m)$
is acyclic. 
Then,

\begin{equation}
0 \longrightarrow K_m \stackrel{\partial_m}{\longrightarrow}
K_{m-1} \longrightarrow \cdots \longrightarrow K_1 \stackrel{\partial_1}{\longrightarrow} K_0 \stackrel{\pi '}{\longrightarrow} R[t] \tag{$\natural$}
\end{equation}

\noindent is also acyclic, where $\pi '$ is the composition of $\pi$ and ${\rm R}(I) \hookrightarrow R[t]$.
Now we take the homogeneous part of $(\natural)$ of degree $m$. 
Then we get an acyclic complex
\[
0 \longrightarrow [K_m]_m \stackrel{[\partial_m]_m}{\longrightarrow}
[K_{m-1}]_m \longrightarrow \cdots \longrightarrow [K_1]_m \stackrel{[\partial_1]_m}{\longrightarrow} [K_0]_m \stackrel{\epsilon_m}{\longrightarrow} R
\]
of finitely generated free $R$-modules, 
where $\epsilon_m$ is the composition of
\[
[\pi ']_m : S_m = [K_0]_m \rightarrow I^mt^m \hspace{1em} and \hspace{1em} I^mt^m \ni at^m \rightarrow a \in R
\]
As a consequence, it follows that $\gr \deti{1}{M} \geq m+1$ by 2.4 and \cite[1.4.13]{bh}.
On the other hand, by taking the homogeneous part of $(\natural)$ of degree 1, we get an acyclic complex
\[
0 \longrightarrow [K_1]_1 \stackrel{[\partial_1]_1}{\longrightarrow} [K_0]_1 \longrightarrow R
\]
of finitely generated free $R$-modules, and so it follows that $\gr \deti{m}{M} \geq 2$ by 2.5. 

In the rest, by induction on $m$, 
we prove $\gr \deti{k}{M} \geq m-k+2$ for $1\leq \forall k \leq m$.
This is certainly true if $m=1$ or 2 by our observation stated above.
So we may assume $m \geq 3$. 
Suppose that we have $\ell := \gr \deti{k}{M} \leq m-k+1$ for some $k$ with
$1 < k < m$. We take a maximal $R$-regular sequence 
$c_1, c_2, \dots, c_{\ell}$ contained in $\deti{k}{M}$.
Then, there exists $\gp \in \ass{R}{R/(c_1, c_2, \cdots, c_{\ell})R}$ 
such that $\deti{k}{M} \subseteq \gp$.
When this is the case, we have $\gr \gp = \ell \leq m-k+1$. 
Hence $\deti{1}{M} \nsubseteq \gp$ as
$\gr \deti{1}{M} \geq m+1 > m-k+1$. 
Then, applying elementary operations to $M$ in $\som{m}{m+1}{R_{\gp}}$,
we get a matrix of the form
\[
\left(\begin{array}{@{\,}c|ccc@{\,}}
1 & 0 & \cdots & 0\\ \hline
0 &&& \\
\vdots & &N &\\
0&&&
\end{array}
\right),
\]
where $N = (y_{ij}) \in \som{m-1}{m}{R_{\gp}}$. 
Let us notice that $\deti{m-1}{N} = \deti{m}{M} R_{\gp} = I_{\gp}$.
As ${\rm S}(I)_{\gp} \stackrel{\sim}{\longrightarrow} {\rm R}(I)_{\gp}$ 
by (i), we have ${\rm S}(I_{\gp}) \stackrel{\sim}{\longrightarrow}
{\rm R}(I_{\gp})$. For $2 \leq \forall i \leq m$, we set
\[
g_i = \displaystyle \sum_{j=1}^{m} y_{i-1,j} T_{j+1} \in S'
 \subseteq R_{\gp}[T_1,T_2, \dots, T_{m+1}] = S_{\gp},
\]
where $S' = R_{\gp}[T_2, \cdots, T_{m+1}]$. 
By 2.3, there exists an isomorphism $\varphi : S_{\gp} \stackrel{\sim}{\to} S_{\gp}$ of
$R_{\gp}$-algebras such that 
$\varphi((f_1, f_2, \dots, f_m) S_{\gp}) = (T_1,g_2, \dots, g_m)S_{\gp}$.

\begin{equation}
\begin{split}
{\rm S}(I_{\gp}) & \cong {\rm S}(I)_{\gp} \\
& \cong S_{\gp}/(f_1, f_2, \dots, f_m)S_{\gp} \,\,\, (\mbox{by (ii)}) \\
& \cong S_{\gp}/(T_1,g_2, \dots, g_m)S_{\gp} \,\,\, (\mbox{isomorphism induced from $\varphi$}) \\
& \cong S'/ (g_2, \dots, g_m)S'. \notag \\
\end{split}
\end{equation}
Moreover, we get $\gr (g_2, \dots, g_m)S' = m-1$. 
Therefore the hypothesis of induction implies that
\[
\gr \deti{k-1}{N} \geq (m-1) - (k-1) + 2 = m - k + 2,
\]
and so $\dep R_{\gp} \geq m-k+2$ as $\deti{k-1}{N} = \deti{k} {M}R_{\gp} \subseteq \gp R_{\gp}$,
which contradicts to $\gr \gp \leq m-k+1$. 
Thus we see $\gr \deti{k}{M} \geq m-k+2$ for $1 \leq \forall k \leq m$
and the proof is complete.

\end{document}